\documentclass{article}
\usepackage[utf8]{inputenc}
\usepackage{amsmath,amssymb,amsthm,bm}
\usepackage{amsfonts}
\usepackage{mathrsfs}
\usepackage{graphicx}
\usepackage{epstopdf}
\usepackage{enumitem}
\usepackage{overpic}
\usepackage{multirow}
\usepackage{nicefrac}
\numberwithin{equation}{section}
\usepackage{subfig}
\usepackage{float}
\usepackage{caption}

\usepackage{multicol}
\usepackage{dsfont}
\usepackage{bbm}
\usepackage{bm}
\usepackage{url}
\usepackage{verbatim}
\usepackage{algorithm,algorithmic}
\usepackage{xcolor}

\usepackage[colorlinks,linkcolor=blue,
            citecolor=blue]{hyperref}
            
\usepackage[
backend=biber,
style=nature,
sorting=ynt
]{biblatex}

\usepackage{geometry}
\usepackage{cleveref}
\usepackage{biblatex}
\usepackage{slashed}

\newcommand{\diag}{\mathrm{diag}}

\newcommand{\rdot}{{\rm rdot}}

\newtheorem{question}{Question}

\providecommand{\keywords}[1]
{
  \small	
  \textbf{\textit{Key words---}} #1
}

\graphicspath{
{images/}
}
\title{Stable Discrete Minimization of  Conformal Energy for Disk Conformal Parameterization}

\author{
Zhong-Heng Tan\thanks{School of Mathematics, Southeast University, Nanjing, 210096, and Nanjing Center for Applied Mathematics, Nanjing, 211135, People’s Republic of China. ({Email: zhtan95math@gmail.com}.)} 
\and 
Zhenyue Zhang\thanks{Corresponding author, Nanjing Center for Applied Mathematics, Nanjing, 211135, and Zhejiang University, Hangzhou, 310027, People’s Republic of China. ({Email: zyzhang@zju.edu.cn})}
}

\date{}

\begin{document}
\maketitle
\begin{abstract}
    Conformal energy minimization is an efficient approach to compute conformal parameterization. In this paper, we develop a stable algorithm to compute conformal parameterization of simply connected open surface, termed Stable Discrete Minimization of  Conformal Energy (SDMCE). The stability of SDMCE is reflected in the guarantee of one-to-one and on-to property of computed parameterization and the insensitivity on the initial value. On one hand, SDMCE can avoid degeneration and overlap of solution, also, SDMCE is folding free. On the other hand, even if given poor initial value, it can still correct it in very little computational time. The numerical experiments indicate SDMCE is stable and competitive with state-of-the-art algorithms in efficiency.
\end{abstract}

\keywords{disk conformal parameterization, conformal energy, adaptive parameters, folding free}

\section{Introduction}

Conformal parameterization aims to find a conformal mapping that maps a given manifold $\mathcal{M}$ to another manifold $\mathcal{M}^*$ with a regular shape, such as square, disk and sphere. 
The conformal mapping, also known as angle-preserving mapping, means the intersection angle of every pair of intersecting arcs on $\mathcal{M}^*$ equals to that of corresponding pre-images on $\mathcal{M}$. 
Conformal parameterization has been widely applied in texture mapping \cite{SHSA00,LBPS02}, image morphing \cite{XFYF15,MHWW17}, medical imaging \cite{XGYW04} and physical analysis \cite{MMKA20}. Therefore, how to compute the conformal mapping efficiently, accurately and stably is a significant problem. 

The computation of discrete conformal mapping has been developed in variant approaches so far. 
Some are to solve equation such as partial differential equations to compute it.
\cite{SHSA00} solves Laplace-Beltrami equation by finite element method to obtain the conformal mapping from closed surface to sphere. 
FLASH \cite{PTKC15}, FDCP \cite{PTLM15} and LDCP \cite{GPLM18} utilize the compound of multiple quasi-conformal mappings to construct a conformal mapping onto sphere/disk. 
Different from these approaches, some are to minimize defined energy funcationals or angle distortions to compute it. 
ABF \cite{SASE01} and ABF++ \cite{SALB05} directly compute the optimal intersection angles by minimizing the angle distortion first and then align the vertices on the planar domain. 
Discrete Ricci flow \cite{MJJK08} and its generalization \cite{YLRG09} compute the conformal mapping by minimizing the Ricci energy, which can be used on surface with arbitrary topologies. 
MIPS \cite{Horm00} minimizes sum of ratio of Dirichlet energy and area on each triangle, called as deformation distortion, which flattens the open surface on the planar domain.
LSCM \cite{LBPS02} and DCP \cite{MDMM02} minimize the conformal energy to compute conformal mapping onto a free boundary domain on $\mathcal{R}^2$. 
Spectral conformal parameterization \cite{PMYT08} transform the conformal energy minimization into a generalized eigenvalue problem. 
Whereas, CEM \cite{MHWW17} and CCEM \cite{YCWW21} compute conformal mapping from simply connected open surface to disk by minimizing the conformal energy. 

In this paper, we focus on the disk conformal parameterization by minimizing the conformal energy, that is, the computation of discrete conformal mapping from open surface to the 2D unit disk. According to the uniformization theorem, the existence of this mapping is guaranteed. Compared with other approaches, conformal energy minimization is efficient since it yields a linear system without extra constraints. However, there are two issues that should be addressed in numerical computation:

\begin{question}
How can we guarantee a computed solution obeying the onto-restriction to the unit disk?
\end{question}

Generally, the boundary vertices are restricted to be unit for mapping the surface to the unit disk. However, this constraint cannot preserve the required onto-property automatically. Actually, minimization of discrete conformal energy under the boundary constraint has a degeneration solution that maps the whole vertices onto a single point with unit length. The centralization strategy used in \cite{} may partially address this degeneration, though this strategy may change the boundary distribution. 

\begin{question}
How can we guarantee a computed solution preserving the one-to-one property? Or equivalently, how can we avoid folding that may appear in the numerical solution?
\end{question}

The folding phenomenon exists in many algorithms \cite{MFKH05,PMYT08}. There are some issues that may result in the foldings.
First, no special strategies are taken into accounted in discretization modeling generally. Second, numerical computation errors may also result in foldings. For instance, it is required to solving a linear system for determining the interior vertices from boundary vertices. However, an ill-conditioned coefficient matrix may yield relative large computational errors. Third, unsuitable distribution of vertices may result in relative large discrete errors even for Delaunay triangulation.

%LSCM and DCP may generate folding triangles \cite{MFKH05}. That is, the one-to-one property cannot be guaranteed. Additionally, spectral conformal parameterization cannot yet guarantee the one-to-one property \cite{PMYT08}, as does CEM and CCEM. 

In this paper, we aim to address these issues. The main technique is to modify discrete conformal energy by adding a penalty term to penalize the deviation of the discrete area to the true one, assuming that the area of the  target surface is known, and the discrete area is computable. It is important to suitably choose the penalty constant for many penalty methods. In our case, we propose an adaptive approach to address the problem of penalty setting. We call the approach as Stable Discrete Minimization of  Conformal Energy (SDMCE) since it can address the onto and folding issues with the three advantages:
\begin{itemize}
    \item SDMCE is insensitive on the initial value and of strong ability to correct for degeneration. Even if given a pool initial value, such as random order boundary vertices and partial degeneration boundary, it has a strong ability on correcting these boundary points.
    \item SDMCE can avoid overlap phenomenon and eliminate both boundary and interior folding triangles. That is, the parameterization computed by SDMCE is one-to-one. 
    \item SDMCE is of fast computation. It has similar computational efficiency to FDCP, LDCP and CCEM, which are very efficient algorithm in recent years. Meanwhile, The computational time does not depend on the initial value.
\end{itemize}

The SDMCE can be easily implemented when the target surface is a 2D unit disk since the true area is known and the discrete area can be measured by the mapping $f$ of the boundary vertices.

\section{Stable discrete minimization of conformal energy}

Theoretically, the continuous conformal energy ${\cal E}_C(f) = {\cal E}_D(f)-{\cal A}^*$ can be approximated by the discrete conformal energy. 
Since the continuous conformal energy is ${\cal E}_C(f) = {\cal E}_D(f)-{\cal A}^*$, we use the discrete conformal energy
\begin{align}\label{def: discrete CE} 
    {\cal E}_C^d(f) = {\cal E}_D^d(f)-{\cal A}^*.
\end{align}
By the convergence analysis for discrete Dirichlet energy, we also conclude that the optimal conformal transformation that minimizes the conformal energy ${\cal E}_C(f)$ can also be converged by the optimal discrete conformal transformation that minimizes the discrete conformal energy
\begin{align}\label{discrete conformal energy}
    \min_{f\in{\cal F}}\Big\{{\cal E}_C^d(f)
    ={\cal E}_D^d(f)-{\cal A}^*\Big\},
\end{align}
or equivalently, minimizes the discrete Dirichlet energy,
as $d\to 0$, provided that the discrete solution $f_d$ guarantees the one-to-one mapping from ${\cal M}$ to ${\cal M}^*$. The discrete conformal energy ${\cal E}_C^d(f)={\cal E}_D^d(f)-{\cal A}^*$ can also be taken as a reliable measurement for the approximation of optimal discrete Dirichlet energy to the optimal continuous Dirichlet energy ${\cal E}_D(f) = {\cal A}^*$.

\subsection{Stable discrete minimization}

Clearly, minimizing the discrete Dirichlet energy or discrete conformal energy is a constrained problem since it asks for a mapping from ${\cal M}$ onto a given ${\cal M}^*$. 
However it is not easy to obey the constraint $f({\cal M}) = {\cal M}^*$ when one looks for a discrete solution. The difficulties involved in numerical computation may include the two issues that are tightly related with each other: \begin{itemize}
    \item the restriction of imaged points $\{f_i = f(v_i)\}$ lying on the surface ${\cal M}^*$, and 
    \item the approximation to the whole surface ${\cal M}^*$ by the piece-wisely linear triangular faces.
\end{itemize} 
The restriction issue may be easily addressed for some simple surfaces such as a 2D disk or a 3D ball. The approximation issue tightly depends on the distribution of $\{f_i\}$ on the surface ${\cal M}^*$. To our acknowledge, no algorithms for discrete conformal transformation can address this issue well. 

The approximation issue is implicitly and partially touched in \cite{MHWW17,YCWW21} by two strategies when ${\cal M}^*$ is a 2D disk: imposing zero mean for $\{f_i\}$ and minimizing a discrete conformal energy $\frac 12 \langle Lf,f\rangle- A(f)$ with a discrete area $A(f)$ in term of $\{f_i\}$ as shown in (\ref{A(f):disk}) later. Starting with a good initial guess that solves a discrete harmonic equation, the approach performs well in some examples. If the initial guess is suitably chosen, a perplexity mistake may happen: a poor solution has a very small value of the measurement. This phenomenon may also happen when the true area ${\cal A}^* = \pi$ is used to replace the estimated area $A(f)$, if only ask $\{f_i\}$ belong to the disk and the boundary points have unit length, together with the zero mean, since the discrete solution may be folded under these restrictions, as shown later in Figure \ref{fig:FaceHo_degen}.

To address the approximation issue, we propose a penalized model of (\ref{discrete conformal energy}). Suppose that the restriction issue on ${\cal M}^*$ is addressed, {\it  i.e.}, $\{f_i\}\subset{\cal M}^*$, and the area ${\cal A}^*$ of ${\cal M}^*$ is known or well estimated.\footnote{These assumptions are easily satisfied for some special ${\cal M}^*$.} 
If no foldings happen on the imaged points $\{f_i\}$, that is, all the triangle vertices $\{f_i,f_j,f_k\}$ have the consistent anticlockwise order,\footnote{Theoretically, the folding should not also happen on the interior points.} the area of ${\cal M}^*$ can be well approximate by the discrete area $A(f)$, {\it i.e.} the area of piecewise constant surface of the triangulation $\{F_{ijk} = [f_i,f_j,f_k]\}$,
\begin{align}\label{area}
    A(f) = \frac{1}{2}\sum_{V_{ijk}\in V}\|(f_i-f_k)\times(f_j-f_k)\|
    =\frac{1}{2}\sum_{V_{ijk}\in V}\sin\theta_{ij}\|f_{jk}\|_2\|f_{ki}\|_2,
\end{align}
where $\theta_{ij}$ is the interior angle opposite to the edgy connection $f_i$ and $f_j$ in the triangle $F_{ijk} = \{f_i,f_j,f_k\}$. Suppose that the approximation of the discrete area $A(f)$ to ${\cal A}^*$ is measured by an error function $\varepsilon_A(f)$. 
We can transform the minimization of discrete conformal energy (\ref{discrete conformal energy}) to the penalized problem \begin{align}\label{discr conformal energy}
    \min_{f_i\in{\cal M}^*} \Big\{\frac 12 \langle Lf,f \rangle -{\cal A}^*
     +\mu \varepsilon_A(f)\Big\}
\end{align}
for stably minimizing the discrete conformal energy, where $f=f(V)$ is the discrete form of $f$ and each row is a point $f_\ell$. We call it as Stable Discrete Minimization of  Conformal Energy (SDMCE). In this paper, we suggest $\varepsilon_A(f) = \frac{1}{2}\big({\cal A}^*- A(f)\big)^2$.

By the way, the true area ${\cal A}^*$ can also be replaced by the estimated discrete ares $A(f)$ if this inequality $A(f)\leq {\cal A}^*$ is always true for $f({\cal M})\subset{\cal M}^*$. That is, (\ref{discr conformal energy}) is approximately equivalent to
\begin{align}\label{discr conformal energy-Af}
    \min_{f_i\in{\cal M}^*} \Big\{\frac 12 \langle Lf,f \rangle - A(f)
     +\mu \varepsilon_A(f)\Big\}.
\end{align}
Generally, 
\begin{align*}
    \min_{f_i\in{\cal M}^*} \Big\{\frac 12 \langle Lf,f \rangle -{\cal A}^*
     +\mu \varepsilon_A(f)\Big\}
     &\leq \min_{f_i\in{\cal M}^*} \Big\{\frac 12 \langle Lf,f \rangle -A(f)
     +\mu \varepsilon_A(f)\Big\}\\
     &\leq \frac 12 \langle LF^*,F^* \rangle -A(f^*)
     +\mu \varepsilon_A(f),
\end{align*}
where $F^* = \{f_i^* = f^*(v_i)\}$ is the optimal solution of (\ref{discr conformal energy}). 

It is a bit complicated to solve (\ref{discr conformal energy}) or (\ref{discr conformal energy-Af}), due to the computational difficulty of the estimated area function $A(f)$. We may have to generate a new approach to estimate the area in terms of the imaged points $\{f_i\}$. This will be one of the topics in our coming work. However, when ${\cal M}^*$ is a 2D disk, this difficulty can be released since a simple formula of $A(f)$ exists. We will propose a modification of SDMCE for the unit 2D disk that releases the onto-restriction in Section \ref{sec:disk}.

\subsection{Adaptive setting of the penalty parameter} 

A suitable value of the penalty parameter $\mu$ should yield a small discrete conformal energy ${\cal E}_C^d(f) = {\cal E}_D^d(f) -{\cal A}^*$ in absolute value, and meanwhile, the area deviation ${\cal A}^*-A(f)$ is also small and positive theoretically. The latter is important to preserve a one-to-one mapping. %Notice that if $A(f)<0$, we should change the signs of $y$-components of $f = [x,y]$, {\it i.e.}, reset $y:=-y$, which changes the ordering of boundary points and guarantees $A(f)>0$, but does not change ${\cal E}_C^d(f)$. 
These two metrics are computable given a solution $f$. Meanwhile, the angle preserving is also checkable via the angle errors $\{\epsilon_{\theta_\ell}(f)\}$ or their average $\bar \epsilon(f)$. Hence, the parameter $\mu$ can be easily tuned to get a solution as good as possible via a simple tuning rule.  Below we show such a tuning role for $\mu$.

\begin{itemize}
    \item Initially, we set $\mu = 0$. If the solution $f$ gives a negative ${\cal E}_C^d(f)$ or a negative $\epsilon_A(f)$, we increase $\mu$ as $\mu:=\mu+s_\mu$ successively, where $s_\mu$ is a constant step length or increased step length step-by-step. We use the later one as $s_\mu:= s_\mu+10$ before $\mu$ is updated. The initial $s_\mu = 0$.
    
    \item As soon as both ${\cal E}_C^d(f)$ and $\epsilon_A(f)$ are nonnegative with a chosen $\mu$, we set $\mu'=\mu$, $f' = f$. Then, modify $s_\mu:=\max\{s_{\mu},10\}$ and increase $\mu'$ to $\mu'' = \mu'+s_\mu$, and compute a solution $f'' = f$ with $\mu=\mu''$. If one of ${\cal E}_C^d(f'')$ and $\epsilon_A(f)$ is negative, $\mu''$ should be also increased as above for $\mu'$.\footnote{In our experiments, this case never happens.}
    
    \item If the gap between $\bar\epsilon(f')$ and $\bar\epsilon(f'')$ is a bit large, we can further tuning the $\mu'$ or $\mu''$ to get a better solution than $f'$ and $f''$. Algorithm \ref{alg:mu} gives the details of the strategy of tuning $\mu$. 
    
    \item By nonnegative ${\cal E}_C^d(f)$ or $\epsilon_A(f)$ numerically, we mean ${\cal E}_C^d(f)>-\tau$ and $\epsilon_A(f)>-\tau$ with a given accuracy $\tau>0$.
\end{itemize}

\begin{algorithm}[t]	
\caption{Adaptive parameter tuning of $\mu$.} \label{alg:mu}
	\begin{algorithmic}[1]
        \REQUIRE Accuracy  $\tau$
        % \ENSURE  
        \STATE Choose initial boundary points with central angles $T_0 = \{t_i\}$ and set $\mu = 0$, $s_\mu = 0$ 
        \WHILE {one of $\mu'$ and $\mu''$ is empty}
        \REPEAT
            \STATE Update $\mu := \mu+s_\mu$ and $s_\mu  = s_\mu+10$, and 
                compute $f$ by SDMCE starting with $T_0$.
            \STATE If $|\epsilon_A(f)|<0.1$, update $T_0$ by the central angles      $T$ of the boundary points of $f$.  
        \UNTIL ${\cal E}_C^d(f)>-\tau$ and $\epsilon_A(f)>-\tau$.
        \STATE If $\mu'$ is empty, set $\mu' = \mu$ and $f' = f$, otherwise
                set $\mu'' = \mu$ and $f''=f$. 
        \ENDWHILE
        \IF{$\bar\epsilon(f')<(1-\tau)\bar\epsilon(f'')$ and $\mu'>0$}
            \WHILE {$s_\mu>5$}
            \REPEAT
                \STATE Update $s_\mu:=s_\mu/2$, $\mu = \lfloor \mu'-s_\mu \rfloor$,
                and compute $f$ by SDMCE again.
            \UNTIL ${\cal E}_C^d(f)>-\tau$ and $\epsilon_A(f)>0$.
            \STATE If $\epsilon_\theta(f)<(1-\tau)\bar\epsilon(f')$, update $\mu'$ and        $f'$ by $\mu$ and $f$ respectively.
            \ENDWHILE
        \ENDIF
        \IF {$\bar\epsilon(f'')<(1-\tau)\bar\epsilon(f')$}
            \REPEAT
            \STATE Update $\mu = \mu''+s_\mu$ and compute $f$ by SDMCE again.
            \STATE If $\epsilon_\theta(f)<(1-\tau)\bar\epsilon(f'')$, 
                update $\mu''$ and $f''$ by $\mu$ and $f$, respectively. 
            \UNTIL{$f''$ does not updated.}
        \ENDIF
        \RETURN $f=f'$ if $\bar\epsilon(f')<\bar\epsilon(f'')$, or $f=f''$ otherwise. 
	\end{algorithmic}
\end{algorithm}

\section{The SDMCE for unit 2D disk}\label{sec:disk}

For the 2D disk ${\cal M}^* = \big\{\mathbf{y}\in {\cal R}^2: \|\mathbf{y}\|_2\leq 1\big\}$, a discrete conformal energy different from (\ref{def: discrete CE})
\begin{align}\label{def:DCE_CCEM}
    \widetilde {\cal E}_{C}^d(f) = \frac 12 \langle Lf,f \rangle - A(f)
\end{align}
was considered in \cite{YCWW21},
where $A(f)$ is the discrete are of the unit disk as (\ref{area}), and it has a simple representation in terms of the boundary points of $\{f_i\}$ as shown in (\ref{A(f):disk}). The approach CCEM proposed in \cite{YCWW21} minimizes $\widetilde {\cal E}_{C}^d(f)$, {\it i.e.}, solves 
\begin{align}\label{CCEM}
    \min \Big\{\frac 12 \langle Lf,f \rangle - A(f)\Big\}
\end{align}
subject to the unit and exactly separate boundary points of $f$ in a correct order. 
%It was shown in \cite{YCWW21} that for a local minimal solution of this constrained problem, the exact separation of boundary points is satisfied. And hence, the separation restriction is released in computation. It is clear that a set of $\{f_i\}$ degenerating to a single point is an optimal solution of (\ref{CCEM}) since the discrete conformal energy defined by (\ref{def:DCE_CCEM}) is zero. To avoiding such a degenerate solution, a good starting set $\{f_i\}$ is required. For instance, starting with the solution of discrete Laplace-Beltrami equation, the CCEM works well on the data sets. 

\subsection{The model of SDMCE for unit disk}

For the unit 2D disk ${\cal M}^* = \big\{\mathbf{y}\in {\cal R}: \|\mathbf{y}\|_2\leq 1\big\}$, the area ${\cal A}^* = \pi$ and its approximate area $A(f)$ is determined by the polygon of boundary points $\{f_i: i\in \Gamma\}$ of $\{f_i\}$. Without loss of generality, we assume that the index set of the boundary points is $\Gamma = \{1,2,\cdots,n\}$ with $n$ is the number of boundary points. 
Representing the boundary points as $f_i = (x_i,y_i)$ with $x_i=\cos t_i$, $y_i = \sin t_i$, a monotone sequence $\{t_i\}\subset[0,2\pi)$ for $i=1,2,\cdots,n$, and $t_0=t_n$, we have 
\begin{align}\label{A(f):disk}
    A(f)=\frac{1}{2}\sum_i\sin(t_{i}-t_{i-1})
    =\frac{1}{2}\sum_i\big(y_ix_{i-1}-x_iy_{i-1}\big)
    =\frac{1}{2}\langle x_\Gamma,D_2y_\Gamma\rangle,
\end{align}
where $D_2$ is a second-order difference operator and skew-symmetric: $D_2^\top = -D_2$, satisfying 
\[
    D_2(y_1,\cdots,y_n)^\top 
    = \big(y_2-y_n,y_3-y_1,\cdots,y_n-y_{n-2},y_1-y_{n-1}\big)^\top.
\]
By $D_2^\top = -D_2$, we also have $A(f)=\frac{1}{2}\langle y_\Gamma,-D_2x_\Gamma\rangle$, and 
\[
    A(f) = \frac{1}{4}\big\langle [x_\Gamma,y_\Gamma],D_2[y_\Gamma,-x_\Gamma]\big\rangle
    = \frac{1}{4}\big\langle f_\Gamma, D_2f_\Gamma \Theta\big\rangle
\]
with $\Theta = \left[{0\atop 1}\ {-1\atop 0}\right]$. Hence, the SDMCE model (\ref{discr conformal energy}) for unit disk becomes
\begin{align}\label{prob:disk}
    \min_{\|f_j\|_2<1,j\notin\Gamma,\|f_i\|_2=1,i\in \Gamma}\Big\{
        \frac12\langle Lf,f\rangle -\pi +\frac{\mu}{2}
        \big(\pi-\frac{1}{4}\langle D_2f_\Gamma\Theta,f_\Gamma\rangle\big)^2\Big\}.
\end{align}

The constrained optimization problem can be solved by classical optimization algorithms such as trust-region method \cite{TCYL96}, gradient descent method \cite{}, or nonlinear conjugated gradient (NCG) \cite{POLYAK69}. We use a modified version of NCG proposed in \cite{GHLT06} in the experiments reported in this paper.  Notice that the objective function has a relative simple structure: the first term is a quadratic form whose coefficient matrix $L$ is diagonal dominant generally if the triangulation is Delaunay \cite{YCWW21}, and the penalization term is the square of a shifted quadratic form. %A iterative scheme via fast matrix iteration can be constructed as follows. 

Let $\Gamma_c$ be the index set of interior points. It is also the complement of $\Gamma$. For simplicity, we release the restriction $\|f_i\|_2<1$ for interior points since it is automatically satisfied, and keep $\|f_i\|_2=1$ for the boundary points in (\ref{prob:disk}). Hence, partitioning the Laplacian matrix $L$ in the $2\times2$  block form with blocks $L_{\Gamma,\Gamma}$, $L_{\Gamma,\Gamma_c}$,  $L_{\Gamma_c,\Gamma}$, and $L_{\Gamma_c,\Gamma_c}$, and $f$ in the two blocks $f_{\Gamma}$ and $f_{\Gamma_c}$, we can rewrite
\[
    \langle Lf,f\rangle
    =\langle L_{\Gamma,\Gamma}f_{\Gamma}+L_{\Gamma,\Gamma_c}f_{\Gamma_c},
        f_{\Gamma}\rangle
    +\langle L_{\Gamma_c,\Gamma}f_{\Gamma}+L_{\Gamma_c,\Gamma_c}f_{\Gamma_c},
        f_{\Gamma_c}\rangle ,
\]
and the restriction on boundary points is $\rdot(f_\Gamma,f_\Gamma) = e$, where $\rdot(A,B)$ for two matrices $A$ and $B$ in the same size is a column vector of the inner productions of corresponding rows of $A$ and $B$, and $e$ is a column vector of all ones in a suitable length.
Let $\eta(f_\Gamma) = \pi-\frac{1}{4}\langle D_2f_\Gamma\Theta,f_\Gamma\rangle$ for simplicity. The KKT conditions of optimal solutions of (\ref{prob:disk}) are  
\begin{align}
    &L_{\Gamma,\Gamma}f_{\Gamma}+L_{\Gamma,\Gamma_c}f_{\Gamma_c}
    -\frac{\mu\eta(f_\Gamma)}{2}D_2f_\Gamma \Theta
    = \Lambda f_\Gamma, \label{KKT:DCE1}\\
    &L_{\Gamma_c,\Gamma}f_\Gamma+L_{\Gamma_c,\Gamma_c}f_{\Gamma_c} = 0,\quad
    \rdot(f_\Gamma,f_\Gamma) = e.\label{KKT:DCE2}
\end{align}
where $\Lambda  = \diag(\lambda)$ with a column vector $\lambda$.
% At a KKT point, 
% $\langle Lf,f\rangle 
%     = \langle L_{\Gamma,\Gamma}f_{\Gamma}+L_{\Gamma,\Gamma_c}f_{\Gamma_c},
%         f_{\Gamma}\rangle$, and the value of objective function is
% \begin{align*}
%     \frac{1}{2}\langle Lf,f\rangle-\pi+\frac{\mu}{2}\eta^2(f_\Gamma)
%     &= \frac{\mu\eta(f_\Gamma)}{4}\langle D_2f_\Gamma \Theta, f_{\Gamma}\rangle
%         +\frac{1}{2}e^\top \lambda -\pi+\frac{\mu}{2}\eta^2(f_\Gamma)\\
%     &=\mu\eta(f_\Gamma)(\pi-\eta(f_\Gamma))
%         +\frac{1}{2}e^\top \lambda -\pi+\frac{\mu}{2}\eta^2(f_\Gamma)\\
%     &= \frac{1}{2}e^\top \lambda -\pi+\mu\eta(f_\Gamma)(\pi-\frac{1}{2}\eta(f_\Gamma)).
% \end{align*}

The restriction 
$L_{\Gamma_c,\Gamma}f_\Gamma+L_{\Gamma_c,\Gamma_c}f_{\Gamma_c} = 0$
for a KKT point shows the dependence of inner points to the boundary points. This restriction is enforced on feasible solution in \cite{YCWW21}. That is, take the boundary points $f_\Gamma$ as variables and $f_{\Gamma_c} =  -L_{\Gamma_c,\Gamma_c}^{-1}L_{\Gamma_c,\Gamma}f_\Gamma$ before optimizing, assuming that $L_{\Gamma_c,\Gamma_c}$ is invertible. In this case,  
\[
    \langle Lf,f\rangle
    =\langle L_{\Gamma,\Gamma}f_{\Gamma}+L_{\Gamma,\Gamma_c}f_{\Gamma_c},
        f_{\Gamma}\rangle
    = \big\langle (L_{\Gamma,\Gamma} - L_{\Gamma,\Gamma_c}L_{\Gamma_c,\Gamma_c}^{-1}L_{\Gamma_c,\Gamma})f_{\Gamma},f_{\Gamma}\big\rangle
    = \langle Sf_{\Gamma},f_{\Gamma}\rangle,
\]
where $S = L_{\Gamma,\Gamma} - L_{\Gamma,\Gamma_c}L_{\Gamma_c,\Gamma_c}^{-1}L_{\Gamma_c,\Gamma}$ is the Schur complement of $L_{\Gamma,\Gamma}$. Hence, the problem (\ref{prob:disk}) becomes 
\begin{align}\label{prob:SDMCE_disk}
    \min_{\|f_i\|_2=1,i\in \Gamma}\Big\{\frac{1}{2}
        \langle Sf_{\Gamma},f_{\Gamma}\rangle-\pi +\frac{\mu}{2}
        \big(\pi-\frac{1}{4}\langle D_2f_\Gamma\Theta, f_\Gamma\rangle\big)^2\Big\}.
\end{align}
% We call it as Stable Minimization of Conformal Energy (SMCE) for unit disk.

The benefit of the above simplified SDMCE for disk is that it significantly reduces the problem scale since the number of boundary points $|\Gamma|$ is significantly smaller than the number of interior points $|\Gamma_c|$. The coefficient matrix $S$ could be constructed explicitly or implicitly for computing $Sf_\Gamma$. The former solves $L_{\Gamma_c,\Gamma_c} H = L_{\Gamma_c,\Gamma}$ and forms $S = L_{\Gamma,\Gamma} -L_{\Gamma,\Gamma_c}H$ once only, and computes $Sf_\Gamma$ in each iteration.
The later solves $L_{\Gamma_c,\Gamma_c} h = L_{\Gamma_c,\Gamma}f_\Gamma$ at first and then $Sf_\Gamma = L_{\Gamma,\Gamma}f_\Gamma - L_{\Gamma,\Gamma_c}h$ in each iteration.
Hence, if an iterative solver used for (\ref{prob:SDMCE_disk}) needs a lot of iterations more than $|\Gamma|$, the explicit strategy costs less than the implicit strategy. In our experiments, we still use the explicit strategy. By the way, computing interior points from boundary points strongly depends on the condition number of $L_{\Gamma_c,\Gamma_c}$; It may loss efficiency when $L_{\Gamma_c,\Gamma_c}$ is ill-conditioned since small perturbation from boundary points will enlarged to the interior points.

\begin{table}[t!]
    \centering
    \resizebox{16cm}{!}{
    \begin{tabular}{|cc|ccc|c|cc|cc|}\hline\hline
        \multicolumn{2}{|c|}{Data} 
            & \multicolumn{3}{c|}{Fixed $\mu = 10$} 
            & \multicolumn{5}{c|}{Variant $\mu\in[a,b]$}\\\hline
        \multirow{2}{*}{ID}& \multirow{2}{*}{Name}
            & \multirow{2}{*}{${\cal E}_C^d(f)$} & \multirow{2}{*}{$\epsilon_\theta(f)$} 
                & \multirow{2}{*}{\!\!Time(s)\!\!} & \multirow{2}{*}{$[a,b]$} 
                & \multicolumn{2}{c|}{${\cal E}_C^d(f)$} & \multicolumn{2}{c|}{$\epsilon_\theta(f)$}\\
                \cline{7-10}
        &&&&&&  Mean     & Std      & Mean     & Std\\\hline
1 & BimbaStatue     &  9.64e-04 & 6.23e-03 &  18.37 & [  0, 500] &  9.64e-04 & 5.26e-09 & 6.23e-03 & 8.32e-08 \\ 
2 & Buddha          &  6.60e-04 & 4.97e-03 & 203.08 & [  0, 250] &  6.60e-04 & 1.38e-12 & 4.97e-03 & 4.58e-10 \\ 
3 & CCH             &  5.42e-03 & 1.16e-02 &  11.51 & [  0, 250] &  5.42e-03 & 1.74e-07 & 1.16e-02 & 9.69e-06 \\ 
4 & CHLin1          &  6.98e-03 & 1.30e-02 &   7.76 & [  0, 250] &  6.98e-03 & 8.70e-06 & 1.31e-02 & 5.68e-05 \\ 
5 & CYHo            &  4.40e-04 & 4.84e-03 &   5.19 & [  0, 250] &  4.40e-04 & 1.39e-11 & 4.84e-03 & 6.95e-09 \\ 
6 & CYHo15          &  5.55e-04 & 5.36e-03 &   4.92 & [  0, 250] &  5.55e-04 & 8.90e-11 & 5.36e-03 & 8.84e-09 \\ 
7 & CYHo22          &  1.22e-02 & 1.84e-02 &   9.41 & [  0, 250] &  1.22e-02 & 4.20e-06 & 1.86e-02 & 1.11e-04 \\ 
8 & CYHo23          &  2.51e-02 & 1.93e-02 &   5.65 & [  0, 250] &  2.51e-02 & 4.95e-07 & 1.93e-02 & 1.44e-05 \\ 
9 & CYHo31          &  2.09e-02 & 1.89e-02 &   6.59 & [  0, 250] &  2.09e-02 & 1.08e-06 & 1.89e-02 & 7.42e-06 \\ 
10 & ChineseLion     &  1.31e-02 & 2.50e-02 &   0.53 & [  0, 250] &  1.31e-02 & 2.50e-07 & 2.50e-02 & 2.20e-06 \\ 
11 & CowboyHat       & -6.80e-03 & 1.47e-02 &   0.03 & [3600,5000] &  1.09e-03 & 6.15e-04 & 3.66e-02 & 1.23e-03\\ 
12 & Dress           &  6.40e-02 & 6.10e-02 &   0.06 & [  0, 250] &  6.42e-02 & 2.44e-04 & 6.09e-02 & 2.04e-05 \\ 
13 & Ear             &  6.99e-02 & 1.04e-01 &   0.01 & [  5,  50] &  7.78e-02 & 1.22e-02 & 1.06e-01 & 1.68e-03 \\ 
14 & Face            &  1.36e-02 & 3.21e-02 &  16.06 & [  0, 500] &  1.36e-02 & 2.67e-05 & 3.21e-02 & 8.96e-06 \\ 
15 & FaceHo          &  1.18e-04 & 2.66e-03 &   9.36 & [  0, 500] &  1.18e-04 & 1.90e-11 & 2.66e-03 & 3.24e-10 \\ 
16 & FaceLin         &  9.27e-04 & 3.31e-03 &  17.93 & [  0, 500] &  9.27e-04 & 3.99e-09 & 3.31e-03 & 5.44e-07 \\ 
17 & Femur           & -2.47e-03 & 1.65e-02 &   0.37 & [ 20, 500] &  2.00e-03 & 7.75e-04 & 1.63e-02 & 6.10e-05 \\ 
18 & Foot            & -1.40e-02 & 2.54e-02 &   0.13 & [150, 200] &  2.18e-04 & 1.49e-04 & 2.52e-02 & 8.51e-06 \\ 
19 & Hand            & -2.18e-03 & 2.12e-02 &   1.76 & [ 15, 500] &  5.27e-03 & 1.96e-03 & 2.09e-02 & 1.48e-04 \\ 
20 & HumanBrain      &  1.26e-03 & 2.66e-02 &   0.90 & [ 11, 500] &  1.43e-02 & 4.45e-03 & 2.60e-02 & 2.80e-04 \\ 
21 & KnitCapMan      &  5.17e-03 & 9.48e-03 &   2.50 & [  0, 500] &  5.17e-03 & 7.78e-07 & 9.48e-03 & 2.06e-06 \\ 
22 & LCH             &  6.98e-03 & 1.30e-02 &   7.50 & [  0, 500] &  6.98e-03 & 1.55e-05 & 1.30e-02 & 1.03e-04 \\ 
23 & LeftHand        & -2.32e-03 & 2.10e-02 &   1.69 & [ 50,4000] &  6.49e-03 & 8.19e-04 & 2.04e-02 & 1.10e-04 \\ 
24 & MaxPlanckD      &  7.22e-03 & 1.05e-02 &   1.06 & [  0, 500] &  7.22e-03 & 3.47e-07 & 1.05e-02 & 3.21e-07 \\ 
25 & Nefertiti       & -2.91e-02 & 5.03e-02 &   0.01 & [110, 200] &  5.68e-03 & 3.20e-03 & 5.03e-02 & 2.16e-03 \\ 
26 & NefertitiStatue &  1.67e-03 & 5.53e-03 & 275.70 & [  0, 500] &  1.67e-03 & 3.09e-09 & 5.53e-03 & 2.84e-08 \\ 
27 & StanfordBunny   &  1.59e-02 & 1.91e-02 &   1.11 & [  0, 500] &  1.59e-02 & 2.19e-06 & 1.91e-02 & 7.17e-07 \\ 
28 & \!\!StanfordBunny2\!\!  
&  8.87e-03 & 1.84e-02 &   1.08 & [  0, 500] &  8.87e-03 & 1.24e-06 & 1.84e-02 & 8.15e-07 \\ 
 \hline\hline
    \end{tabular}}
    \caption{Efficiency of the SDMCE starting with an equal distance boundary points on 28 real world data sets, measured by the discrete energy ${\cal E}_C^d(f)$, the average angle error $\epsilon_\theta(f)$, and computational time}
    \label{tab:stablity}
\end{table}

%\begin{figure}[t]
%    %\centering
%    \hspace{-50pt}\includegraphics[width=1.2\linewidth,height = 2.5cm]{images/Ear.eps}
%    
%    \hspace{-50pt}\includegraphics[width=1.2\linewidth,height = 2.5cm]{images/Hand.eps}
%    
%    \hspace{-50pt}\includegraphics[width=1.2\linewidth,height = 2.5cm]{images/Nefertiti.eps}
%    \caption{Dependence of the solutions of SDMCE on the penalty parameter $\mu$, measured by the discrete conformal energy ${\cal E}_C^d(f)$ (left), area deviation $\pi-A(f)$ (middle), and average angle error $\bar \epsilon(f)$ (right) on {\em Ear, Hand}, and {\em Nefertiti}, starting with equal-distantly distributed boundary points on the unit circle.}
%    \label{fig:negative energy}
%\end{figure}

The SDMCE model (\ref{discr conformal energy}) or (\ref{prob:SDMCE_disk}) for unit disk works very well. It can efficiently avoid the degeneration of solution. For instance, starting with a set of equal-distance points on the unit circle, if we do not use the penalty term, or set $\mu=0$ equivalently, the minimized discrete Dirichlet energy may be wrongly zero on some data sets
%the six data sets {\em Ear, Femur, Foot, Hand, HumanBrain, LeftHand} and {\em Nefertiti} 
-- the solutions degenerate to a single point together. However, starting the same initial points, the SDMCE with a not very small $\mu$ gives a good solution on these data sets as shown in the section of numerical experiments.     

\begin{table}[t!]
    \centering
    \resizebox{16cm}{!}{
    \begin{tabular}{|c|c@{\, }c|c@{\, }c|c@{\, }c|c@{\, }c|c@{\, }c|c@{\, }c|} \hline\hline
    \multirow{3}{*}{Data} & \multicolumn{10}{c|}{True ordering with variant ratio $\rho = 2\pi/\ell$} 
        & \multicolumn{2}{c|}{Wrong ordering}\\\cline{2-13}
        & \multicolumn{2}{c|}{$\rho=0.4$} & \multicolumn{2}{c|}{$\rho=0.8$}
        & \multicolumn{2}{c|}{$\rho=1.2$} & \multicolumn{2}{c|}{$\rho=1.6$}
        & \multicolumn{2}{c|}{$\rho=2.0$} 
        & \multirow{2}{*}{${\cal E}_C^d(f)$} & \multirow{2}{*}{$\epsilon_\theta(f)$}\\\cline{2-11}
        & ${\cal E}_C^d(f)$ & $\epsilon_\theta(f)$ & ${\cal E}_C^d(f)$ & $\epsilon_\theta(f)$ 
        & ${\cal E}_C^d(f)$ & $\epsilon_\theta(f)$ & ${\cal E}_C^d(f)$ & $\epsilon_\theta(f)$ 
        & ${\cal E}_C^d(f)$ & $\epsilon_\theta(f)$ &&\\\hline
BimbaStatue     &  9.58e-04 & 6.23e-03 &  1.09e-03 & 6.29e-03 &  1.05e-03 & 6.26e-03 &  9.51e-04 & 6.23e-03 &  9.56e-04 & 6.24e-03 &  9.47e-04 & 6.24e-03\\
Buddha          &  1.65e-03 & 5.32e-03 &  7.32e-04 & 5.01e-03 &  6.82e-04 & 4.98e-03 &  8.65e-04 & 5.03e-03 &  6.97e-04 & 4.99e-03 &  6.74e-04 & 4.98e-03\\
CCH             &  1.10e-02 & 1.21e-02 &  5.42e-03 & 1.16e-02 &  5.43e-03 & 1.15e-02 &  5.41e-03 & 1.16e-02 &  5.42e-03 & 1.17e-02 &  1.10e-02 & 1.20e-02\\
CHLin1         &  9.86e-03 & 1.39e-02 &  6.97e-03 & 1.30e-02 &  6.94e-03 & 1.29e-02 &  6.95e-03 & 1.30e-02 &  9.86e-03 & 1.37e-02 &  9.86e-03 & 1.40e-02\\
CYHo            &  3.97e-04 & 4.83e-03 &  4.17e-04 & 4.81e-03 &  4.33e-04 & 5.10e-03 &  5.56e-04 & 5.83e-03 &  3.99e-04 & 4.95e-03 &  9.39e-04 & 7.21e-03\\
CYHo15         &  4.77e-04 & 5.54e-03 &  6.51e-04 & 5.58e-03 &  5.51e-04 & 5.58e-03 &  6.62e-04 & 6.30e-03 &  4.82e-04 & 5.52e-03 &  8.19e-04 & 7.55e-03\\
CYHo22         &  1.22e-02 & 1.87e-02 &  1.22e-02 & 1.87e-02 &  2.20e-02 & 1.78e-02 &  1.83e-02 & 1.83e-02 &  1.22e-02 & 1.86e-02 &  1.83e-02 & 1.82e-02\\
CYHo23         &  1.51e-02 & 2.06e-02 &  1.51e-02 & 2.05e-02 &  1.89e-02 & 1.97e-02 &  1.89e-02 & 1.98e-02 &  1.51e-02 & 2.05e-02 &  1.51e-02 & 2.05e-02\\
CYHo31         &  1.90e-02 & 1.97e-02 &  1.90e-02 & 1.97e-02 &  2.10e-02 & 1.88e-02 &  2.10e-02 & 1.87e-02 &  1.90e-02 & 1.97e-02 &  4.06e-03 & 2.33e-02\\
ChineseLion     &  1.31e-02 & 2.50e-02 &  1.31e-02 & 2.50e-02 &  1.31e-02 & 2.50e-02 &  1.31e-02 & 2.50e-02 &  1.31e-02 & 2.50e-02 &  1.31e-02 & 2.50e-02\\
CowboyHat       &  5.75e-05 & 3.45e-02 &  5.44e-05 & 3.45e-02 &  1.17e-01 & 4.39e-02 &  1.17e-01 & 4.31e-02 &  1.18e-01 & 4.23e-02 &  1.41e-01 & 5.36e-02\\
Dress           &  6.40e-02 & 6.10e-02 &  6.40e-02 & 6.10e-02 &  6.40e-02 & 6.10e-02 &  6.40e-02 & 6.10e-02 &  6.40e-02 & 6.10e-02 &  6.40e-02 & 6.10e-02\\
Ear             &  6.98e-02 & 1.04e-01 &  6.99e-02 & 1.04e-01 &  6.99e-02 & 1.04e-01 &  6.99e-02 & 1.04e-01 &  6.99e-02 & 1.04e-01 &  6.99e-02 & 1.04e-01\\
Face            &  1.36e-02 & 3.21e-02 &  1.36e-02 & 3.20e-02 &  1.35e-02 & 3.21e-02 &  1.36e-02 & 3.21e-02 &  1.36e-02 & 3.21e-02 &  1.41e-02 & 3.26e-02\\
FaceHo          &  6.33e-05 & 4.35e-03 &  9.63e-05 & 2.90e-03 &  1.66e-04 & 2.78e-03 &  3.50e-04 & 3.32e-03 &  1.50e-04 & 4.05e-03 &  1.24e-04 & 3.52e-03\\
FaceLin         &  1.05e-03 & 3.43e-03 &  1.06e-03 & 3.47e-03 &  9.88e-04 & 3.45e-03 &  1.07e-03 & 3.59e-03 &  1.03e-03 & 3.35e-03 &  1.04e-03 & 3.67e-03\\
Femur           &  3.12e-04 & 1.64e-02 &  2.59e-04 & 1.64e-02 &  2.51e-04 & 1.64e-02 &  3.12e-04 & 1.64e-02 &  3.82e-04 & 1.64e-02 &  2.46e-04 & 1.64e-02\\
Foot            &  2.60e-04 & 2.52e-02 &  2.56e-04 & 2.52e-02 &  2.63e-04 & 2.52e-02 &  2.59e-04 & 2.52e-02 &  2.64e-04 & 2.52e-02 &  2.59e-04 & 2.52e-02\\
Hand            &  5.12e-04 & 2.11e-02 &  5.23e-04 & 2.11e-02 &  5.00e-04 & 2.11e-02 &  4.88e-04 & 2.11e-02 &  4.06e-04 & 2.11e-02 &  4.62e-04 & 2.11e-02\\
HumanBrain      &  1.10e-03 & 2.66e-02 &  1.29e-03 & 2.66e-02 &  1.26e-03 & 2.66e-02 &  1.17e-03 & 2.66e-02 &  1.19e-03 & 2.66e-02 &  1.12e-03 & 2.66e-02\\
KnitCapMan      &  5.17e-03 & 9.48e-03 &  5.17e-03 & 9.49e-03 &  5.17e-03 & 9.49e-03 &  5.17e-03 & 9.49e-03 &  5.17e-03 & 9.48e-03 &  5.17e-03 & 9.49e-03\\
LCH             &  9.86e-03 & 1.39e-02 &  6.97e-03 & 1.30e-02 &  6.94e-03 & 1.29e-02 &  6.95e-03 & 1.30e-02 &  9.86e-03 & 1.37e-02 &  6.95e-03 & 1.32e-02\\
LeftHand        &  4.19e-03 & 2.07e-02 &  4.34e-03 & 2.07e-02 &  4.23e-03 & 2.07e-02 &  4.32e-03 & 2.07e-02 &  4.26e-03 & 2.07e-02 &  4.24e-03 & 2.07e-02\\
MaxPlanckD      &  7.21e-03 & 1.05e-02 &  7.21e-03 & 1.05e-02 &  7.22e-03 & 1.05e-02 &  7.21e-03 & 1.05e-02 &  7.21e-03 & 1.05e-02 &  7.21e-03 & 1.05e-02\\
Nefertiti       &  1.32e-03 & 4.74e-02 &  1.32e-03 & 4.74e-02 &  1.31e-03 & 4.74e-02 &  1.32e-03 & 4.74e-02 &  1.32e-03 & 4.74e-02 &  1.32e-03 & 4.74e-02\\
NefertitiStatue &  1.68e-03 & 5.53e-03 &  1.82e-03 & 5.54e-03 &  1.74e-03 & 5.53e-03 &  1.72e-03 & 5.53e-03 &  1.68e-03 & 5.53e-03 &  1.69e-03 & 5.53e-03\\
StanfordBunny   &  1.59e-02 & 1.91e-02 &  1.59e-02 & 1.91e-02 &  1.59e-02 & 1.91e-02 &  1.59e-02 & 1.91e-02 &  1.59e-02 & 1.91e-02 &  1.59e-02 & 1.91e-02\\
StanfordBunny2  &  8.87e-03 & 1.84e-02 &  8.87e-03 & 1.84e-02 &  8.87e-03 & 1.84e-02 &  8.87e-03 & 1.84e-02 &  8.87e-03 & 1.84e-02 &  8.87e-03 & 1.84e-02\\\hline\hline
    \end{tabular}}
    \caption{The stability of SDMCE on the wrong starting points.}
    \label{tab:wrong initial points}
\end{table}

\subsection{Advantages of SDMCE for disk}

\begin{figure}[t]
    % \centering
    \hspace{-50pt}\includegraphics[width=1.2\linewidth,height = 5cm]{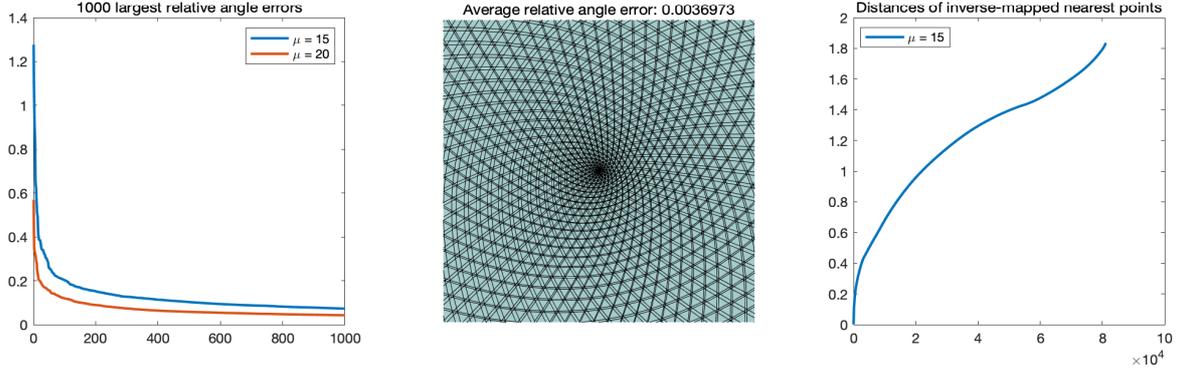}
    \caption{A phenomenon of degenerated solution on {\em FaceHo}: the disk is covered twice by the solution of SDMCE with $\mu=15$, starting with boundary points along the unit circle twice. The degeneration disappears when $\mu\geq 17$. Left: distribution of 1000 largest relative angle errors. Middle: partial domain of the solution with $\mu=15$. Right: the distribution of distances $\{d_i\}$ for checking the one-to-one. }
    \label{fig:FaceHo_degen}
\end{figure}

The SDMCE model has some advantages for finding a discrete conformal mapping as shown below. 
\begin{enumerate}
    \item {\bf Fast computation}. It can be solved by NCG or other classical optimization problem as the algorithm CCEM given in \cite{YCWW21}. In our experiments on 28 real world data sets, starting with equal-distantly distributed boundary points on the unit circle for each example, the NCG can give an good solution with computational time slightly less than that of CCEM, due to the penalty to the area deviation of solutions. 
%Table \ref{tab:stablity} shows the discrete conformal energy ${\cal E}_C^d(f)$ defined as (\ref{def:discrete CE_disk}), the mean of relative angle errors $\epsilon_\theta(f)$, and the computational time for each of the 28 data sets, using the same $\mu = 10$. 
    
    \item {\bf Insensitivity on the parameter setting}. The penalty model is not sensitive to the penalization parameter $\mu$ variant in a large range if initially using the equal-distantly distributed boundary points. 
    %Table \ref{tab:stablity} also lists the mean values and the standard deviation of ${\cal E}_C^d(f)$ and $\epsilon_\theta(f)$, among the solutions with the parameter $\mu$ varying in the shown range. These ranges could be much larger than the tested ones.    
    
    \item {\bf Ability of avoiding negative discrete conformal energy}. The computed solution may give a negative discrete conformal energy ${\cal E}_C^d(f)$ as shown in Table \ref{tab:stablity}. Negative ${\cal E}_C^d(f)$ implies a partially degenerate solution with a deficient area covering. Increasing the penalty parameter $\mu$ can decrease the deficiency. Figure \ref{fig:negative energy} illustrates the improvement on {\em Ear, Hand}, and {\em Nefertiti} when $\mu$ is increased. The improvement on {\em Femur, Foot, LeftHand}, and {\em HumanBrain} is similar and omitted. As we mentioned before, minimizing the discrete conformal energy, {\it i.e.}, setting $\mu=0$ in SDMCE, will give a solution degenerated in a single point on these seven examples. 
    Slightly increasing $\mu$ can avoid the degeneration efficiently. As $\mu$ increases, the negative discrete conformal energy becomes to positive, and the covering of target domain is improved significantly. Notice that it also improves the angle preserving of solutions. 
    
    \item {\bf Global convergence}. The penalty model can adaptively correct a wrong initial setting. To show the advantage, we test two kinds of initial setting for boundary points, each of which results in a very poor initial guess. 
    
    \begin{itemize}
        \item[(a)] The initial setting  keeps a correct neighboring order but the ordering points cover an arc length $\ell$ much less than or over larger than $2\pi$. That is, we choose the boundary points $f_i = (\cos t_i,\sin t_i)$ with the ordered central angles $\{t_i\}$ distributed in the interval $[0,\ell]$ with equal neighbor gaps. 
 
        \item[(b)] The boundary points are uniformly distributed in the unit circle with a random order. The random ordering means a seriously heavy folding of these points. 
    \end{itemize}
    
    The SDMCE has a strong ability on correcting these wrong starting boundary points. It can avoid the degeneration on all the tested data sets. Table \ref{tab:wrong initial points} lists the discrete conformal energy, the mean of angle errors with the starting points. By the way, the computational time does not depend the starting points.

    \item {\bf One-to-one property}. A small average value of relative angle errors does not always means a good solution since the one-to-one property between ${\cal M}$ and ${\cal M}^*$ may be lost. The one-to-one can be checked by the distribution of distances $\{d_i\}$, defined at each point $f_i$ of a solution by 
    \begin{align}
        d_i = \|v_i-v_{i'}\|, 
    \end{align}
    where $v_i$ and $v_{i'}$ are the original points of $f_i$ and $f_{i'}$, respectively, and $f_{i'}$ is the nearest one of $f_i$. 
    Figure \ref{fig:FaceHo_degen} illustrates the degeneration happened on {\em FaceHo} in which the target domain is covered twice by the solution (see the middle panel), starting with the boundary points set as (a) with $\rho=2\pi/\ell = 2$ and $\mu\leq16$ is used. The required one-to-one property does not preserved in this case, thought the average angle error is as small as the average error of a good solution. This kind of degeneration disappears when we increase the penalty parameter to $\mu=17$.
\end{enumerate}

The degeneration illustrated by the middle penal of Figure \ref{fig:FaceHo_degen} will also result in a positive discrete conformal energy ${\cal E}_C^d(f)\approx\pi$ since the discrete Dirichlet energy is computed twice, {\it i.e.}, ${\cal E}_D^d(f)\approx 2\pi$. Meanwhile, the area deviation $\epsilon_A(f) = \pi-A(f)\approx -\pi$. Similar phenomenon happens for multiple times of repeated covering. Hence, it is easy to check this kind of degeneration using ${\cal E}_C^d(f)$ or $\epsilon_A(f)$.

The algorithm for adaptively tuning $\mu$ works well and is also stable for starting boundary points. Notice that there are not big additional costs for turning $\mu$. \begin{itemize}
    \item It is not required to modify the Laplacian matrix $L$ or the Schur complement $S$ that is the main cost in solving (\ref{prob:SDMCE_disk}).
    
    \item As soon as $\mu$ is updated, it is required to apply the NCG again. Since (\ref{prob:SDMCE_disk}) is very stable on $\mu$, the solution corresponding to the previous $\mu$ is a very good starting point for reapplying the NCG. That means, the NCG convergs quickly within few of iterations. 
    
    \item The iteration of turning $\mu$ terminates quickly on most of the tested data sets. Only on {\em CowboyHat}, the turning number is larger than that on other data sets, but the additional cost is ignorable. 
\end{itemize}

Table \ref{tab:adaptive mu} shows the results of SDMCE with adaptively chosen $\mu$, when the starting points are equally distributed with correct ordering  or randomly chosen with completely wrong ordering.   

\begin{table}[t!]
    \centering
    \begin{tabular}{|c|c@{\, }c@{\, }c@{\, }c@{\, }c|c@{\, }c@{\, }c@{\, }c@{\, }c|}\hline\hline
    \multirow{2}{*}{Data} 
        & \multicolumn{5}{c|}{Equally distributed initial points}
        & \multicolumn{5}{c|}{Randomly chosen initial points}\\\cline{2-11}
        & $\mu$ & ${\cal E}_C^d(f)$ & $\epsilon_A(f)$ & $\epsilon_\theta(f)$ &\!\!Time(s)\!\! 
        & $\mu$ & ${\cal E}_C^d(f)$ & $\epsilon_A(f)$ & $\epsilon_\theta(f)$ & \!\!Time(s)\!\! \\\hline
BimbaStatue     &    0 & 9.64e-04 & 6.40e-05 & 6.23e-03 &  21.73 &   20 & 9.50e-04 & 6.73e-05 & 6.23e-03 &  24.57\\
Buddha          &    0 & 6.60e-04 & 2.67e-05 & 4.97e-03 & 118.16 &   20 & 6.68e-04 & 3.04e-05 & 4.98e-03 & 160.60\\
CCH             &    0 & 5.42e-03 & 3.23e-04 & 1.16e-02 &  10.34 &   20 & 5.41e-03 & 2.96e-04 & 1.16e-02 &  10.99\\
CHLin1         &   10 & 6.98e-03 & 3.03e-04 & 1.30e-02 &   7.84 &   20 & 6.95e-03 & 2.80e-04 & 1.30e-02 &  13.13\\
CYHo            &   10 & 4.40e-04 & 5.06e-05 & 4.84e-03 &   5.84 &  270 & 1.99e-05 & 1.68e-03 & 6.79e-03 &  14.24\\
CYHo15         &   10 & 5.55e-04 & 4.57e-05 & 5.36e-03 &   5.99 &  270 & 2.13e-04 & 1.56e-03 & 7.21e-03 &  12.66\\
CYHo22         &   10 & 1.22e-02 & 9.20e-04 & 1.87e-02 &  10.33 &   20 & 1.22e-02 & 9.55e-04 & 1.86e-02 &   7.08\\
CYHo23         &    0 & 2.51e-02 & 7.21e-05 & 1.93e-02 &   6.43 &   10 & 1.89e-02 & 4.16e-04 & 1.97e-02 &   7.60\\
CYHo31         &    0 & 2.09e-02 & 7.30e-05 & 1.90e-02 &   7.20 &   90 & 1.28e-02 & 5.50e-03 & 2.17e-02 &  14.78\\
ChineseLion     &   10 & 1.31e-02 & 7.57e-04 & 2.50e-02 &   0.56 &   20 & 1.31e-02 & 7.67e-04 & 2.50e-02 &   0.81\\
CowboyHat       & 3577 & 1.83e-05 & 9.13e-03 & 3.44e-02 &   0.08 & 3568 & 1.26e-05 & 9.13e-03 & 3.44e-02 &   0.14\\
Dress           &   10 & 6.40e-02 & 1.16e-02 & 6.10e-02 &   0.06 &   20 & 6.40e-02 & 1.16e-02 & 6.10e-02 &   0.06\\
Ear             &   10 & 6.99e-02 & 6.78e-02 & 1.04e-01 &   0.02 &   10 & 6.98e-02 & 6.79e-02 & 1.04e-01 &   0.02\\
Face            &    0 & 1.36e-02 & 9.31e-05 & 3.21e-02 &  19.28 &   20 & 1.36e-02 & 1.01e-04 & 3.21e-02 &  19.38\\
FaceHo          &   10 & 1.18e-04 & 2.03e-05 & 2.66e-03 &  10.93 &   20 & 1.31e-04 & 1.39e-04 & 4.07e-03 &  24.02\\
FaceLin         &    0 & 9.27e-04 & 1.98e-05 & 3.31e-03 &  20.53 &  180 & 1.30e-03 & 2.43e-03 & 5.70e-03 &  40.16\\
Femur           &   70 & 2.07e-03 & 4.61e-03 & 1.63e-02 &   0.58 &   70 & 2.07e-03 & 4.59e-03 & 1.63e-02 &   0.53\\
Foot            &  210 & 5.22e-04 & 1.11e-02 & 2.52e-02 &   0.21 &  210 & 5.23e-04 & 1.11e-02 & 2.52e-02 &   0.20\\
Hand            &  150 & 6.34e-03 & 3.10e-03 & 2.08e-02 &   3.26 &  110 & 5.99e-03 & 3.90e-03 & 2.08e-02 &   2.79\\
HumanBrain      &   40 & 1.69e-02 & 3.57e-03 & 2.59e-02 &   1.35 &   40 & 1.69e-02 & 3.59e-03 & 2.59e-02 &   1.20\\
KnitCapMan      &   10 & 5.17e-03 & 4.84e-04 & 9.48e-03 &   2.43 &   50 & 5.17e-03 & 4.88e-04 & 9.48e-03 &   2.94\\
LCH             &   10 & 6.98e-03 & 3.03e-04 & 1.30e-02 &   8.15 &   10 & 6.94e-03 & 2.91e-04 & 1.30e-02 &   9.57\\
LeftHand        &  130 & 5.50e-03 & 3.17e-03 & 2.06e-02 &   2.87 &   70 & 4.81e-03 & 5.04e-03 & 2.07e-02 &   2.69\\
MaxPlanckD      &   10 & 7.22e-03 & 6.74e-04 & 1.05e-02 &   1.19 &   10 & 7.21e-03 & 6.79e-04 & 1.05e-02 &   1.34\\
Nefertiti       &  100 & 2.92e-04 & 2.71e-02 & 4.67e-02 &   0.03 &  100 & 2.92e-04 & 2.71e-02 & 4.67e-02 &   0.04\\
NefertitiStatue &    0 & 1.67e-03 & 4.99e-05 & 5.53e-03 & 330.47 &   10 & 1.68e-03 & 4.96e-05 & 5.53e-03 & 402.01\\
StanfordBunny   &   10 & 1.59e-02 & 8.37e-04 & 1.91e-02 &   3.24 &   20 & 1.59e-02 & 8.36e-04 & 1.91e-02 &   1.43\\
StanfordBunny2  &    0 & 8.87e-03 & 6.82e-04 & 1.84e-02 &   1.73 &   20 & 8.87e-03 & 6.95e-04 & 1.84e-02 &   5.20\\
\hline\hline
    \end{tabular}
    \caption{Comparison of SDMCE with adaptively chosen $\mu$, starting with boundary points equally distributed in the correct ordering (left half) or randomly chosen in a wrong ordering (right half).}
    \label{tab:adaptive mu}
\end{table}

\section{Folding free}

Given a solution $\{f_\ell\}$, folding happens on the solution if there is at least a triangle $F_{ijk}$ of $\{f_\ell\}$ contains a vertex $f_\ell$ differ from the vertices of $F_{ijk}$. It contains the special case when a triangle folds over a boundary edge. That is, the mapping $f_\ell$ of an interior vertex $v_\ell$ appears out of the boundary formed by the boundary points $f_\Gamma$. %The stable model (\ref{prob:SDMCE_disk}) cannot avoid the folding, 

\begin{figure}[t]
    \centering
    \begin{tabular}{c c}
         \includegraphics[width = 7cm]{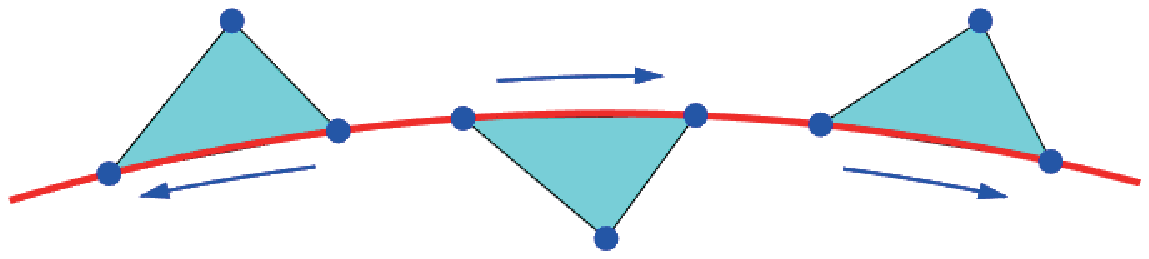} & 
         \includegraphics[width = 7cm]{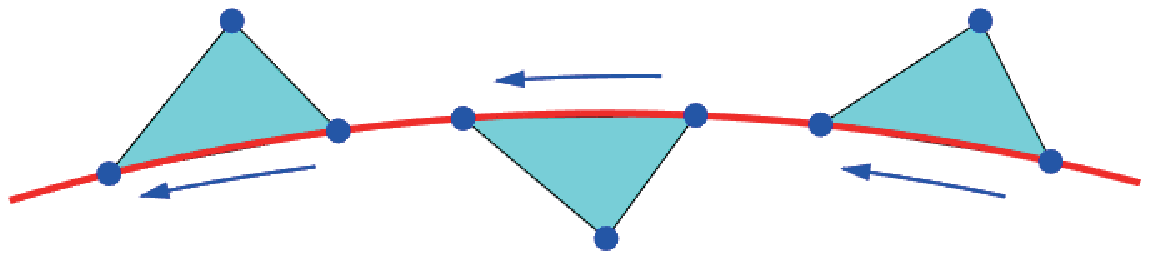}\\
         (a) & (b)
    \end{tabular}
    \caption{Efficiency illustration of boundary penalty. (a): Three kinds of folded boundary triangles without penalty. (b): The first class (left) still exists, the second one (middle) disappears, and the third one (right) becomes to the first class when the penalty strategy is used. } 
    \label{fig:Folding_BT}
\end{figure}

We consider two kinds of folding for the unit disk: boundary folding and triangle folding. By boundary folding, we mean that the central angle $t_{i+1}$ of the boundary vertex $f_{j_{i+1}}$ is smaller than the central angle $t_i$ of $f_{j_i}$ within mod $2\pi$. The triangle folding includes interior triangle folding and boundary triangle folding. By interior triangle folding, we mean that an interior triangle $F_{ijk}$ has its vertices $f_i,f_j,f_k$ in clockwise, that is, its algebra area $\det\big([f_{ji}, f_{kj}]\big)$ is negative. The boundary triangle folding is a bit complicated -- it contains three cases: 
\begin{itemize}
    \item[(1)] A boundary triangle folds over its boundary edge, and the boundary vertices are not folded. Hence, the algebra area is positive.
    \item[(2)] A boundary triangle does not fold over its boundary edge, but the boundary vertices are folded. In this case, the algebra area is also positive.
    \item[(3)] A boundary triangle folds over its folded boundary edge -- the algebra area is positive.
\end{itemize} 
The panel (a) of Figure \ref{fig:Folding_BT} illustrates the three kinds of folding of boundary triangles. Hence, when we account the number of folding triangles according to the negativity of algebra areas, the boundary triangles folding over folded boundary edges are lost. 
%For distinguishing the folded boundary triangles, we call them as negative triangle folding or positive triangle folding, depending on the sign of algebra area.  

The boundary folding can be addressed by optimization. In the next subsection, we consider a modified version of (\ref{prob:SDMCE_disk}) for addressing this kind of folding by adding an adaptive penalty on the boundary points. 
However, the triangle folding is hard to address via optimization methods. 
Generally, the triangle folding mainly due to the ill-conditional coefficient matrix $A$ -- relative large computational error in $S$ or the inverse of $L_{\Gamma_c,\Gamma_c}$ will result in perturbation of the interior points $f_{\Gamma_c}$. After the next subsection, we will also show how to handle the interior folding technically in the later part of this subsection when the boundary folding has been addressed. 

\subsection{Folding free for boundary vertices}

Let $f_{\Gamma} = \{f_{j_i}\}$, where $f_{j_i} = (\cos t_i, \sin t_i)$ with the central angles $\{t_i\}$. We assume that original boundary vertices $\{v_{j_i}\}$ are in anticlockwise order. If no folding happens on $f_\Gamma$, each sector area with arc length of adjacent points $f_{j_i}$ and $f_{j_{i-1}}$, say
\[
    \frac{1}{2}\sin(t_i-t_{i-1}) 
    = \frac{1}{2}\big(\sin t_i\cos t_{i-1}-\cos t_i\sin t_{i-1}\big)
    = f_{j_i}\Theta f_{j_{i-1}}^\top
\]
must be nonnegative, where $j_0 = j_{|\Gamma|}$. Hence, a negative $f_{j_i}\Theta f_{j_{i-1}}^\top$ means the folding on boundary vertices. We add the penalty term 
\[
    \sum_{i=1}^{|\Gamma|} \alpha_i \max\{-f_{j_i}\Theta f_{j_{i-1}}^\top,0\}
    =\alpha^\top \max\big\{\rdot(D_1f_\Gamma\Theta,f_\Gamma),0\big\}
    =\alpha^\top \rdot(D_1f_\Gamma\Theta,f_\Gamma)_+
\]
into (\ref{prob:SDMCE_disk}), where $D_1 = [e_2,\dots,e_{|\Gamma|},e_1]$. It yields the following regularization problem
\begin{align}\label{prob:SDMCE_diskF}
    \min_{\|f_i\|_2=1,i\in \Gamma}\Big\{\frac{1}{2}
        \langle Sf_{\Gamma},f_{\Gamma}\rangle-\pi +\frac{\mu}{2}
        \big(\pi-\frac{1}{4}\langle D_2f_\Gamma\Theta, f_\Gamma\rangle\big)^2
        +\alpha^\top \rdot(D_1f_\Gamma\Theta,f_\Gamma)_+\Big\}.
\end{align}
The NCG can also work well on solving the above problem. Since the main cost is the  computation of the Laplacian matrix $L$ and the Schur complement $S$, the computational cost for solving (\ref{prob:SDMCE_diskF}) is similar as that for (\ref{prob:SDMCE_disk}).

The penalty vector $\alpha$ is nonnegative and adaptively set during the iteration of  the algorithm solving this problem in the rule: Starting with $\alpha=0$ with all zeros, if $f_{j_i}\Theta f_{j_{i-1}}^\top$ is negative, we modify $\alpha_i:=\alpha_i+\delta$, where $\delta$ is a small positive constant, say $\delta = \frac{|\Gamma|}{|\Gamma_c|}$. The penalty action $\alpha_i$ disappears if $f_{j_i}\Theta f_{j_{i-1}}^\top$ becomes to nonnegative. 

Taking this strategy of setting $\alpha$, the penalty method (\ref{prob:SDMCE_diskF}) performs very well on the five data sets on which the boundary folding happens. The folding of boundary vertices disappears on all these data set. That is, the second kind of folding boundary triangles disappears, and the third one becomes to the first one. See the panel (b) of Figure \ref{fig:Folding_BT} for the improvement. Unfortunately, the interior triangle folding and/or the negative boundary triangle folding still exist. The first two lines of Table \ref{tab:folding T} compare the efficiency of the SDMCE using adaptive turning of $\mu$ on {\em Face} when the penalty strategy is used or not. The total 44 folded boundary vertices are corrected, but there are also folded triangles (23 interior triangles and 76 boundary triangles in the first kind).

\begin{table}[t]
    \centering
    \begin{tabular}{|c|ccc|@{\, }c@{\, }|c@{\, }c@{\, }c@{\, }c@{\, }c|}\hline\hline
 \multirow{2}{*}{Strategy} & \multirow{2}{*}{${\cal E}_C^d(f)$} 
    & \multirow{2}{*}{$\epsilon_A(f)$} & \multirow{2}{*}{$\epsilon_\theta(f)$} 
    & Folded &\multicolumn{5}{c|}{Folded triangles in}\\
    & & & & bound. v & total & interior  & 1st b. & 2nd b. & 3th b.\\\hline
 No & 1.36e-02 & 9.3116e-05 & 3.2052e-02 &  44 & 120 &  24 & 52 & 21 & 23\\
 Boundary V & 1.36e-02 & 9.2971e-05 & 3.2052e-02 &   0 &  99 &  23 & 76 &  0 &  0\\
 Boundary T & 1.36e-02 & 9.2971e-05 & 3.2644e-02 &   0 &   8 &   8 &  0 &  0 &  0\\
 Interior T & 1.36e-02 & 9.2971e-05 & 3.2703e-02 &   0 &   0 &   0 &  0 &  0 &  0\\
\hline\hline
    \end{tabular}
    \caption{Efficiency of adaptive SDMCE with the strategies for folding free on {\em Face}}
    \label{tab:folding T}
\end{table}

\subsection{Folding free for boundary triangles}

As soon as the boundary folding is address, there is only the first kind of folding of boundary triangles, though the number of this folded boundary triangles is increased. In this subsection, we show how to address the boundary triangle folding.  

Let $T_{ijk} = \{v_i,v_j,v_k\}$ be a boundary triangle of ${\cal M}$ with boundary edgy $e_{jk} = (v_j,v_k)$, corresponding to the folded boundary triangle of a given solution $\{f_\ell\}$. Let $V_i = \{v_{i_\ell}\}$ be the adjacent neighboring vertices including the boundary vertices $v_j$ and $v_k$. We estimate the weights $\{w_{i\ell}\}$ of $v_i$ in an approximate convex combination in terms of its adjacent neighbors by solving 
\begin{align}\label{prob:convW}
    \min_{\{w_{i\ell}\}}
    \Big\|v_i - \sum_{\ell:v_\ell\in V_i}w_{i\ell}v_{i_\ell}\Big\|_2,\quad
    {\rm s.t.}\quad w_{i\ell}\geq0,\ \sum_{\ell:v_\ell\in V_i}w_{i\ell}=1.
\end{align}
It is not difficult to solve this problem. As soon as the solution $w_i = \{w_{i\ell}\}$ is available, $f_i$ is updated by 
\begin{align}\label{update f_i}
    f_i: = \sum_{\ell:v_\ell\in V_i}w_{i\ell}f_{i_\ell}.
\end{align}

The updating rule (\ref{update f_i}) based on (\ref{prob:convW}) for all first kind of folded boundary triangles can easily be addressed. Due to the correction on this kind of boundary triangles, some folded interior triangles connected with the corrected vertices are also corrected simultaneously.
For instance, the 76 folded boundary triangles in the first kind are corrected, and meanwhile, there are 15 interior triangles are also unfolded. This is an interesting advantage of the above approach for unfolding boundary triangles.

\subsection{Folding free for interior triangles}

The problem of unfording interior triangles is a bit complicated, compared with that for boundary triangles, since (1) it is not clear which vertex of the folded interior triangle results in the folding, and (2) folded vertices of some folded interior triangles may be connected with each others. 

Let $T_{ijk} = \{v_i,v_j,v_k\}$ be a interior triangle of ${\cal M}$ corresponding to a folded interior triangle of a given solution $\{f_\ell\}$. Let $V_i' = \{v_{i_\ell}\}$ $V_j' = \{v_{j_\ell}\}$ and $V_k' = \{v_{k_\ell}\}$ be the adjacent neighboring vertices, not including the vertices $v_i$, $v_j$, and $v_k$. As in the above subsection, we can obtain the weight vectors $w_i = (w_i', w_{ij},w_{ik})$, $w_j = (w_j', w_{jk},w_{ji})$, and $w_k = (w_k', w_{ki},w_{kj})$, each solves a similar minimization problem as (\ref{prob:convW}) for representing $v_i$, $v_j$, and $v_k$, in terms of their connected neighbors. Hence, we have  
\[
    \left[\begin{array}{ccc}
         1 & -w_{ij} & -w_{ik}\\
         -w_{ji} & 1 & -w_{jk}\\
         -w_{ki} & -w_{kj} & 1
    \end{array}\right]
    \left(\begin{array}{c}
    v_i \\ v_j \\ v_k
    \end{array}\right)\approx 
    \left[\begin{array}{c}
    w_i'V_i'\\ w_j'V_j' \\ w_k'V_k'
    \end{array}\right].
\]
We unfold the folded interior triangle $F_{ijk} = \{f_i,f_j.f_k\}$ by updating it to
\begin{align}\label{update F_{ijk}}
    \left(\begin{array}{c}
    f_i \\ f_j \\ f_k
    \end{array}\right)
    = \left[\begin{array}{ccc}
         1 & -w_{ij} & -w_{ik}\\
         -w_{ji} & 1 & -w_{jk}\\
         -w_{ki} & -w_{kj} & 1
    \end{array}\right]^{-1}
    \left[\begin{array}{c}
    w_i'F_i'\\ w_j'F_j' \\ w_k'F_k'
    \end{array}\right],
\end{align}
where $F_i'$, $F_j'$, and $F_k'$ are the connected neighbors of $F_{ijk}$ as $V_i'$, $V_j'$, and $V_k'$.

\begin{algorithm}[t]
\caption{SDMCE for disk parameterization of open surfaces}
\begin{algorithmic}[1]
    \REQUIRE vertices $\{v_i\}$, indices $T = \{T_{ijk}\}$ of triangles, initial 
    \STATE Determine index set $\Gamma$ of boundary vertices, initially set $f_\Gamma$, $\alpha_i = 0$, $c = |\Gamma|/|T|$.
    \STATE Construct the Laplace matrix $L$ and the Schur complement $S$.
    \STATE Apply Algorithm \ref{alg:mu} to determine $\mu$ and a solution $\{f_i\}$. 
    \WHILE {boundary vertex folding happens}
        \STATE Update $\alpha_i := \alpha_i+c$ for folded boundary vertices, 
        and solve (\ref{prob:SDMCE_diskF}) starting with the current $\{f_i\}$. 
    \ENDWHILE
    \WHILE {boundary triangle folding happens}
        \STATE Solve (\ref{prob:convW}) and update each interior vertex of the triangles.
    \ENDWHILE
    \WHILE {interior triangle folding happens}
        \STATE Compute the weight vectors $w_i,w_j,w_k$ by solving (\ref{prob:convW}) and update each folded triangle as (\ref{update F_{ijk}}).
    \ENDWHILE
\end{algorithmic}
\end{algorithm}

Applying the updating rule on all the folded interior triangles, these triangles are unfolded generally. Since the modification on the vertices of the folded interior triangles may slightly lead to folding of connected triangles, or partial of the interior triangles are not completely unfolded via one step of the unfolding, we might repeat this approach on the remaining folded triangles. For instance, this approach can unfold all the 8 interior triangles, however, three newly folded interior triangles appear. Applying (\ref{update F_{ijk}}) on these three interior triangles, no other folding appears. The last row of Table \ref{tab:folding T} shows the result when the (\ref{update F_{ijk}}) is applied twice on {\em Face}.

\section{Numerical experiments and comparisons}

In this section we will show the numerical behaviors of the algorithm SDMCE, and compare it with algorithms FDCP \cite{PTLM15}, LDCP \cite{GPLM18}, and CCEM \cite{YCWW21} on 28 real-world data sets listed in Table \ref{tab:}. The SMDCE model is solved via the modified NCG given in \cite{}. We use the relative angle change of each interior angle $\{\theta_\ell^v\}$ of triangles
\begin{align}\label{error:angle}
    \epsilon_{\theta_\ell}(f) = \frac{|\theta_{\ell}^v-\theta_{\ell}^f|}{\theta_{\ell}^v}
\end{align}
to measure the angel preserving, where each angle $\theta_{\ell}^v$ is measured in the anticlockwise order of the original triangle $V_{ijk} = \{v_i,v_j,v_k\}$, and so does $\theta_{\ell}^f$ of the transformed triangle $F_{ijk} = \{f_i,f_j,f_k\}$. The mean of all the relative angle errors is denoted as $\epsilon_{\theta}(f)$. We will also use the discrete conformal energy on the unit disk
\begin{align}\label{def:discrete CE_disk} 
    {\cal E}_C^d(f) = {\cal E}_D^d(f)-{\cal A}^* = \frac{1}{2}\langle Lf,f\rangle -\pi
\end{align}
to measure the approximation of discrete conformal energy of the solution to the minimal conformal energy for the ideal conformal transformation. Clearly, the error (\ref{def:discrete CE_disk}) can also be used to measure the approximation error of discrete Dirichlent energy ${\cal E}_D^d(f)$ to the continuous Dirichlet energy ${\cal E}_D(f) = {\cal A}^*$ for the ideal conformal transformation,
\[
    {\cal E}_D^d(f) - {\cal E}_D(f) = {\cal E}_D^d(f)-{\cal A}^*.
\]
Here we do not take the absolute value for the approximation error since as we show later, a negative error ${\cal E}_C^d(f)$ implies a solution not good enough for preserving angles.
We do not use an estimate area of ${\cal A}^*$ in a discretion of conformal energy such as the commonly used error ${\cal E}_D^d(f)-A(f)$ in the literature since such an estimate area $A(f)$ depends on a solution and may result in a small ${\cal E}_D^d(f)-A(f)$ for a degenerate solution $f$. 

\subsection{Efficiency of the penalty model}

\subsection{Stability of SDMCE}

We first show the performance of the SDMCE on angel preserving

At first, the 
In our experiments on 28 real world data sets, starting with equal-distantly distributed boundary points on the unit circle for each example, the NCG can give an good solution within a time slightly less than the time of CCEM, due to the penalty to the area deviation of solutions. Table \ref{tab:stablity} shows the discrete conformal energy ${\cal E}_C^d(f)$ defined as (\ref{def:discrete CE_disk}), the mean of relative angle errors $\epsilon_\theta(f)$, and the computational time for each of the 28 data sets, using the same $\mu = 10$.

\subsection{Efficiency  of the folding-free technique} 

\subsection{Comparisons}

In this subsection we show the comparisons of the SDMCE with the algorithms FDCP \cite{PTLM15}, LDCP \cite{GPLM18}, and CCEM \cite{YCWW21} on all the 28 data sets. Below we briefly describe these three  algorithms.

The CCEM solves (\ref{CCEM}) subject to $\|f_i\|^2=1$ for boundary vertices by the quasi-Newton method. It was shown in \cite{YCWW21} that for a local minimal solution of this constrained problem, the exact separation of boundary points is satisfied. And hence, the separation restriction is released in computation. It is clear that a set of $\{f_i\}$ degenerating to a single point is an optimal solution of (\ref{CCEM}) since the discrete conformal energy defined by (\ref{def:DCE_CCEM}) is zero. To avoiding degeneration, a good initial solution is required. 

Different from the CCEM for a conformal mapping directly, both the FDCP and LDCF look for a conformal $f = g_2\circ g_1$ via two {\it quasi-conformal} mappings $g_1: {\cal M}\to D$ and $g_2:D\to D$, and boundary-to-boundary, such that $g_1^{-1}$ and $g$ have equal Beltrami coefficients \cite[Theorem 1]{PTKC15}, {\it i.e.}, both $g_1^{-1}$ and $g_2$ satisfy the Beltrami equation
\begin{align}\label{Beltrami eq}
    \frac{\partial g_1^{-1}}{\partial \bar z}
    = b(z) \frac{\partial g_1^{-1}}{\partial z},\quad
    \frac{\partial g_2}{\partial \bar z}
    = b(z) \frac{\partial g_2}{\partial z},\quad %\forall z\in D,\quad
    g_2\big|_{\partial D} = \partial D,
\end{align}
%for a Beltrami coefficient function $b(z)$, 
where the 2D variable $(x,y)$ is taken as a complex variable $z = x+iy$. Since the $g_1$ could be arbitrarily chosen, the key is to solve the Beltrami equation for $g_2$ given the Beltrami coefficient $b(z)= \frac{\partial g_1^{-1}}{\partial \bar z}/\frac{\partial g_1^{-1}}{\partial z}$ of $g_1^{-1}$. However, the nonlinear boundary restriction $g_2\big|_{\partial D} = \partial D$ increases the difficulty of solving.

The FDCP and LDCP take different strategies to address the difficulty. Basically, FDCP iteratively determines $g_2$ given $g_1$. It takes Cayley transformation $c(z) = i\frac{1+z}{1-z}$ to linearize the boundary restriction, which keeps the Beltrami equation unchanged, and uses the reflection $z\to 1/\bar z$, centralization, and normalization for correcting the solution nearby the pole point $(1,0)$ of Cayley transformation. LDCP looks for an approximately conformal $g_1$ and modifies the boundary condition to $g_2\big|_{\partial D} = g_1\big|_{\partial {\cal M}}$ for simplicity. Here $g_1$ is obtained by applying the couple-quasi-conformal method for conformally mapping a genus-0 closed surface to the unit sphere ${\cal S}^2$ \cite{PTKC15}, copying the open surface to form a required genus-0 closed surface. The similar difficulty from the restriction is address via the stereographic projection, and using M\"obius transformation and normalization for correction. Though the normalization on the boundary points can obey the strict boundary restriction, it may result in folding of boundary triangles. 

By the way, the CCEM solves a inhomogeneous Laplace-Beltrami equation with a right-hand function for an initial solution. Both FDCP and LDCP solve a homogeneous Laplace-Beltrami equation -- with a restriction on boundary points whose distribution of arc length angles equal to the boundary vertices for the open surface in FDCP. The SDMCE simply uses an initial solution with equal-distance distributed arc length angles of boundary points. 

\begin{figure}
    \hspace{-50pt}
    \includegraphics[width = 19cm]{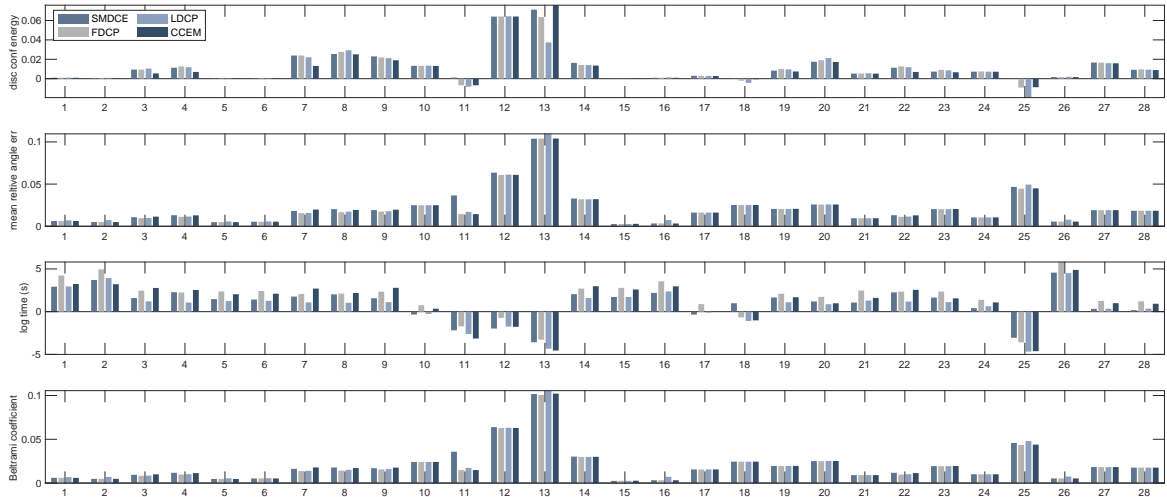}
    \caption{Comparisons of the algorithms SDMCE, FDCP, LDCP, and CCEM on the 28 data sets. From top to bottom: discrete conformal energy, %error of discrete area, 
    average value of relative angle errors, computation time in log, and average value of Beltrami coefficient.}
    \label{fig:comparison}
\end{figure}

Figure \ref{fig:comparison} shows the comparisons of the algorithms SDMCE, FDCP, LDCP, and CCEM on discrete conformal energy, %error of discrete area,
average value of relative angle errors, and computational time of the computed solution on the 28 data sets. 
The SDMCE provides a competitive result on each of the tested data sets. The FDCP, LDCP, and CCEM give negative discrete conformal energies on {\em CowboyHat, Foot, Nefertiti} with the identity number 11, 18, and 25, respectively. The SDMCE costs much time than the LDCP and CCEM on the two largest data sets {\em Buddha} and {\em NefertitiStatue} since it uses the explicit version of the Schur matrix $S$. Using the implicit form of $S$, the SDMCE costs the same as the CCEM. 

\begin{figure}[h]
    \hspace{-50pt}
    \includegraphics[width = 19cm]{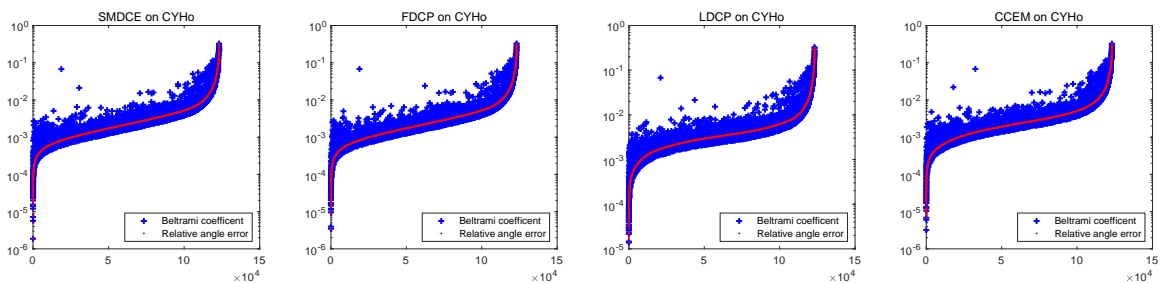}
    \caption{Sorted relative angle errors and the corresponding Distribution of Beltrami coefficient function of the four solutions on the data set {\em CYHo}.}
    \label{fig:BC-AE}
\end{figure}

We also compare the mean of Beltrami coefficient of solutions since the Beltrami coefficient of a conformal mapping should be zero. The FDCP or LDCP could be taken as an iterative method for minimizing the Beltrami coefficient of the compound function $g_2\circ g_1$. It is interesting that the mean of Beltrami coefficient is approximately equal to the mean of relative angle errors, though they are different. See Figure \ref{fig:BC-AE} for the difference on solutions of the four algorithms, where the Beltrami coefficient function is plotted corresponding to the sorted relative angle errors. 

\section{Conclusions}

\printbibliography

@article{GPLM18,
Author = {Gary Pui-Tung Choi and Lok Ming Lui},
Title = {A linear formulation for disk conformal parameterization of simply-connected open surfaces},
Journal = {Advances in Computational Mathematics},
Year = {2018},
Volume = {44},
Number = {1},
Pages = {87--114},
DOI = {10.1007/s10444-017-9536-x},
}

@InProceedings{Horm00,
  author    = {Kai Hormann},
  title     = {MIPS: An efficient global parametrization method},
  booktitle = {ACM Press/Addison-Wesley Publishing Co},
  year      = {2000},
}

@article{LBPS02,
  author     = {L{\'e}vy Bruno and Petitjean Sylvain and Ray Nicolas and Maillot J{\'e}rome},
  title      = {Least Squares Conformal Maps for Automatic Texture Atlas Generation},
  journal    = {ACM Transactions on Graphics},
  year       = {2002},
  volume     = {21},
  number     = {3},
  pages      = {362--371},
  doi        = {10.1145/566654.566590},
}

@article{MDMM02,
  Author = {Mathieu Desbrun and Mark Meyer and Pierre Alliez},
  Title = {Intrinsic parameterizations of surface meshes},
  Journal = {Computer Graphics Forum},
  Year = {2002},
  Volume = {21},
  Number = {3},
  Pages = {209--218},
  DOI = {10.1111/1467-8659.00580},
}

@InProceedings{MFKH05,
  author    = {Michael S. Floater and Kai Hormann},
  title     = {Surface Parameterization: a Tutorial and Survey},
  booktitle = {Advances in Multiresolution for Geometric Modelling},
  year      = {2005},
  pages     = {157--186},
  publisher = {Springer Berlin Heidelberg},
  doi       = {10.1007/3-540-26808-1_9},
}

@article{MHWW17,
  author = {Mei-Heng Yueh and Wen-Wei Lin and Chin-Tien Wu and Shing-Tung Yau},
  title = {An Efficient Energy Minimization for Conformal Parameterizations},
  journal = {Journal of Scientific Computing},
  volume = {73},
  number = {1},
  pages = {203--227},
  year = {2017},
  doi = {10.1007/s10915-017-0414-y},
}

@article{MJJK08,
Author = {Miao Jin and Junho Kim and Feng Luo and Xianfeng Gu},
Title = {Discrete surface Ricci flow},
Journal = {IEEE Transactions on Visualization and Computer Graphics},
Year = {2008},
Volume = {14},
Number = {5},
Pages = {1030--1043},
DOI = {10.1109/TVCG.2008.57},
}

@article{PMYT08,
  author  = {Patrick Mullen and Yiying Tong and Pierre Alliez and Mathieu Desbrun},
  title   = {Spectral conformal parameterization},
  journal = {Computer Graphics Forum},
  year    = {2008},
  volume  = {27},
  number  = {5},
  pages   = {1487--1494},
}

@article{PTKC15,
author = {Pui Tung Choi and Ka Chun Lam and Lok Ming Lui},
title = {FLASH: Fast Landmark Aligned Spherical Harmonic Parameterization for Genus-0 Closed Brain Surfaces},
journal = {SIAM Journal on Imaging Sciences},
volume = {8},
number = {1},
pages = {67-94},
year = {2015},
doi = {10.1137/130950008},
}

@article{PTLM15,
Author = {Pui Tung Choi and Lok Ming Lui},
Title = {Fast Disk Conformal Parameterization of Simply-Connected Open Surfaces},
Journal = {Journal of Scientific Computing},
Year = {2015},
Volume = {65},
Number = {3},
Pages = {1065--1090},
DOI = {10.1007/s10915-015-9998-2},
}

@article{SALB05,
Author = {Alla Sheffer and Bruno L{\'e}vy and Maxim Mogilnitsky and Alexander Bogomyakov},
Title = {ABF++: Fast and robust angle based flattening},
Journal = {ACM Transactions on Graphics},
Year = {2005},
Volume = {24},
Number = {2},
Pages = {311--330},
DOI = {10.1145/1061347.1061354},
}

@article{SASE01,
Author = {Alla Sheffer and Eric de Sturler},
Title = {Parameterization of faceted surfaces for meshing using angle-based
   flattening},
Journal = {Engineering with Computers},
Year = {2001},
Volume = {17},
Number = {3},
Pages = {326--337},
DOI = {10.1007/PL00013391},
}

@article{SHSA00,
Author = {Steven Haker and Sigurd Angenent and Allen Tannenbaum and Ron Kikinis and Guillermo Sapiro and Michael Halle},
Title = {Conformal surface parameterization for texture mapping},
Journal = {IEEE Transactions on Visualization and Computer Graphics},
Year = {2000},
Volume = {6},
Number = {2},
Pages = {181--189},
DOI = {10.1109/2945.856998},
}

@article{XFYF15,
  doi = {10.1007/s00371-015-1188-6},
  year = {2015},
  volume = {32},
  number = {9},
  pages = {1191--1203},
  author = {Xin Fan and Yuyao Feng and Zhi Chai and Xianfeng David Gu and Zhongxuan Luo},
  title = {Image morphing with conformal welding},
  journal = {The Visual Computer}
}

@article{XGYW04,
Author = {Xianfeng Gu and Yalin Wang and Tony F. Chan and Paul M. Thompson and Shing-Tung Yau},
Title = {Genus zero surface conformal mapping and its application to brain surface mapping},
Journal = {IEEE Transactions on Medical Imaging},
Year = {2004},
Volume = {23},
Number = {8},
Pages = {949--958},
DOI = {10.1109/TMI.2004.831226},
}

@article{YCWW21,
  author = {Yueh-Cheng Kuo and Wen-Wei Lin and Mei-Heng Yueh and Shing-Tung Yau},
  title = {Convergent Conformal Energy Minimization for the Computation of Disk Parameterizations},
  journal = {SIAM Journal on Imaging Sciences},
  volume = {14},
  number = {4},
  pages = {1790--1815},
  year = {2021},
  doi = {10.1137/21M1415443},
}

@article{YLRG09,
Author = {Yong-Liang Yang and Ren Guo and Feng Luo and Shi-Min Hu and Xianfeng Gu},
Title = {Generalized Discrete Ricci Flow},
Journal = {Computer Graphics Forum},
Year = {2009},
Volume = {28},
Number = {7},
Pages = {2005--2014},
DOI = {10.1111/j.1467-8659.2009.01579.x},
}

@ARTICLE{MMKA20,  
author={Mirazimi, Maedeh Sadat and Kiyoumarsi, Arash}, 
journal={IEEE Transactions on Transportation Electrification}, 
title={Magnetic Field Analysis of SynRel and PMASynRel Machines With Hyperbolic Flux Barriers Using Conformal Mapping},  
year={2020},  
volume={6},  
number={1},  
pages={52-61},  
doi={10.1109/TTE.2019.2959400}
}

@article{TCYL96,
author = {Thomas F. Coleman and Yuying Li},
title = {An Interior Trust Region Approach for Nonlinear Minimization Subject to Bounds},
journal = {SIAM Journal on Optimization},
volume = {6},
number = {2},
pages = {418-445},
year = {1996},
doi = {10.1137/0806023}
}

@article{POLYAK69,
title = {The conjugate gradient method in extremal problems},
journal = {USSR Computational Mathematics and Mathematical Physics},
volume = {9},
number = {4},
pages = {94-112},
year = {1969},
issn = {0041-5553},
doi = {https://doi.org/10.1016/0041-5553(69)90035-4},
author = {B.T. Polyak}
}

@article{GHLT06,
author = {Gaohang Yu and L\:utai Guan},
title = {Modified PRP methods with sufficient descent property and their convergence properties},
journal = {Acta Scientiarum Naturalium Universitatis Sunyatseni (Chinese)},
volume = {45},
number = {4},
pages = {11-14},
year = {2006}
}
\end{document}